\documentclass[12pt]{amsart}
\usepackage{amsmath}
\usepackage{amssymb}
\usepackage{amsthm}
\usepackage{epsfig,graphicx}
\usepackage{hyperref}
\usepackage{cleveref}

\newtheorem{lem}{Lemma}[section]

\newtheorem{thm}{Theorem}
\newtheorem{corollary}[thm]{Corollary}

\newtheorem{coro}[lem]{Corollary}
\newtheorem{prop}[lem]{Proposition}

\theoremstyle{definition}

\newtheorem{defn}[lem]{Definition}

\newtheorem{exa}[lem]{Example}

\theoremstyle{remark}
\newtheorem{rmk}[lem]{Remark}

\renewcommand{\setminus}{\smallsetminus}

\newcommand{\Pp}{{\mathbb P}}

\newcommand{\nd}{\#_{\operatorname{disks}}}
\newcommand{\ddd}{\Delta^{\operatorname{disks}}}

\newcommand{\cp}{{\mathbb C}{\mathbb P}}

\newcommand{\qg}{{\mathbb Q}\Gamma}
\newcommand{\Q}{{\mathbb Q}}
\newcommand{\R}{{\mathbb R}}
\newcommand{\Z}{{\mathbb Z}}
\newcommand{\C}{{\mathbb C}}
\newcommand{\xx}{\mathcal X}
\newcommand{\dd}{\partial}
\newcommand{\mr}{M_{\R}}
\newcommand{\mc}{M\otimes\C^\times}
\newcommand{\nr}{N_{\R}}
\newcommand{\ff}{\mathcal F}
\newcommand{\Int}{\operatorname{Int}}
\newcommand{\Log}{\operatorname{Log}}

\begin{document}
\title{Singular symplectic spaces and holomorphic membranes}
\author{Sergey Galkin and Grigory Mikhalkin}
\address{PUC-Rio, Matem\'atica, rua Marqu\^es de S\~ao Vicente 225, G\'avea, Rio de Janeiro, Brasil}
\address{HSE University, Russian Federation}
\email{singsimp@galkin.org.ru}
\address{Universit\'e de Gen\`eve,  Math\'ematiques, rue du Conseil-G\'en\'eral 7-9, 1205 Gen\`eve, Suisse}
\email{grigory.mikhalkin@unige.ch} 
\thanks{Research is supported in part by 
the SNSF-grants 200400, 204125
and CNPq-grant PQ 315747. The paper was written during the research visit of the second-named author to the Euler International Mathematical Institute in St. Petersburg, partially supported by a grant from the Government of the Russian Federation, agreements 075-15-2019-1620 and 075-15-2022-289.
}

\begin{abstract}
We set up a topological framework for degenerations of symplectic manifolds into singular spaces paying a special attention to the behavior of Lagrangian manifolds and their (holomorphic) membranes. We show that degenerations into singular toric varieties provide a source of exotic Lagrangian tori.
\end{abstract}

\maketitle

\section{Introduction}
By their very definition, symplectic manifolds are smooth manifolds with a distinguished closed non-degenerate 2-form. Important class of examples is provided by complex projective varieties. If an algebraic subvariety $X^n\subset\cp^N$ is non-singular then it is a $(2n)$-dimensional smooth manifold from the point of view of differential topology. The ambient projective space $\cp^{N}$ comes with the Fubini-Study 2-form $\omega$. As $\omega$ is closed an non-degenerate, $\cp^N$ is a $(2N)$-dimensional symplectic manifold, while $X$ is a $(2n)$-dimensional symplectic manifold inside.

This construction gives an abundant source of examples of symplectic manifolds. But as is, it excludes singular projective manifolds that may be also very useful for constructions in classical symplectic geometry. For instance, suppose that we have a Lagrangian submanifold $L$ inside the smooth locus $X^\circ$ of a singular projective variety $X$. If $X$ is smoothable to a non-singular variety $X_t$ (inside the ambient space $\cp^N$) then $X^\circ$ symplectically embeds to $X_t$ and thus $L$ gets embedded as a Lagrangian variety also to $X_t$. Thus a Lagrangian manifold inside a singular projective varieties yields a Lagrangian in a conventional symplectic manifold after smoothing. 

In their treatment of mirror symmetry for minuscule varieties \cite{BoGa},
Bondal and Galkin suggested an approach via Lagrangian tori obtained from toric degenerations as above.
In 2010 Galkin and Usnich \cite{GaUs} have suggested to use this approach for smoothings of weighted projective planes $\Pp(a^2,b^2,c^2)$ given by the Markov triples $a,b,c\in \Z_{>0}$ with $a^2+b^2+c^2=3abc$. As it was shown by Hacking and Prokhorov \cite{HaPr}, all these planes admit $\Q$-Gorenstein smoothings to the standard plane $\cp^2=\Pp(1,1,1)$. As all other weighted projective planes, the planes $\Pp(a^2,b^2,c^2)$ contain the unique monotone Lagrangian fiber torus (called the {\em Clifford torus}). Under the smoothing of $\Pp(a^2,b^2,c^2)$ to $\cp^2$, it yields a monotone torus $L(a^2,b^2,c^2)\subset\cp^2$.
It was conjectured in \cite{GaUs} that this approach gives non-symplectomorphic  pairs $(\cp^2,L(a^2,b^2,c^2))$ for different Markov triples $(a,b,c)$ as the corresponding Maslov index 2 disk counting is different\footnote{The paper \cite{GaUs} also contains the answer for this disk counting, conjectured at the time of its writing.}.
In particular, $\cp^2$ admit infinitely many different Lagrangian embeddings of tori. This conjecture was verified in \cite{Vianna1}, but in a way circumnavigating the direct algebraic smoothings of $\Pp(a^2,b^2,c^2)$. 
The main theorem of the present paper implies that the tori $L(a^2,b^2,c^2)$ are symplectically different directly from the fact that they come as $\Q$-Gorenstein smoothings of the Clifford tori in toric surfaces $\Pp(a^2,b^2,c^2)$ (given by different toric fans).

The theorem is very robust and general, it describes the shape (Newton polytope) of a Floer potential for monotone Lagrangian tori obtained from any $\Q$-Gorenstein degeneration of Fano manifolds
(of any dimension)
in terms of the combinatorial data of the singular toric variety: Newton polytope of the potential equals to ``fan polytope'' of the degenerate variety.
In the very particular case that singular fiber admits a small resolution the theorem reduces to a theorem of Nishinou--Nohara--Ueda \cite{NNU2,NNU}
used by \cite{BoGa}. 

While the small resolution degenerations are quite rare, just finitely many in a given dimension,
none in dimension two,
the degenerations covered by Theorem \ref{thm-main} form a much larger class. 
For example, a conjecture of Belmans--Galkin--Mukhopadhyay \cite{BGM20}
states that the toric varieties associated with quantum Clebsch--Gordan polytope $\Delta_\gamma$
(given by spherical triangle inequalities and parity conditions)
admit a small resolution if and only if the trivalent graph $\gamma$ is $3$-connected,
so Nishinou--Nohara--Ueda's theorem is applicable only in the case of $3$-connected graphs and is conditional to BGM small resolution conjecture.
On the other hand, our method is applicable for any $\gamma$, see \cite{BGM-sympl}.

Our result also proves a conjecture of Batyrev \cite{Ba04,BK13} in a generalized form: for any $\Q$-Gorenstein degeneration
of a Fano manifold to a toric Fano variety with terminal singularities, the Laurent polynomial prescribed by Batyrev
is mirror dual to the initial manifold.

Heuristics saying that Newton polytope of a Laurent polynomial mirror dual to a Fano manifold
equals to the fan polytope of its toric degeneration appear in \cite{Ga-phd,Ga-ECM,Ga-sigma},
and our result gives a (Floer-)theoretical foundation for them.

In Section 2 we introduce the notion of singular symplectic spaces and its special case of $\qg$-space, suitable for consideration of the canonical class with $\Q$-coefficients. Singular projective varieties provide examples of singular symplectic spaces while $\Q$-Gorenstein projective varieties provide examples of $\qg$-spaces. In the same section we consider the corresponding notions of smoothing and $\qg$-smoothing. In Section 3 we specialize these notions to the toric set-up while in section 4 we study holomorphic membranes of Maslov index 2 in a $\qg$-smoothing of a normal toric variety $X_\Delta$ and prove that they determine the convex hull of the primitive vectors in the toric fan of $X_\Delta$. 

In the special case of toric degenerations of $\cp^2$ to $\Pp(a^2,b^2,c^2)$ treated in \cite{HaPr}, this implies that pairs $(\cp^2,L(a^2,b^2,c^2))$ for different Markov triple $(a,b,c)$ are non-symplectomorphic, as predicted in \cite{GaUs}.
Similarly, our result can be applied to 
Manon's toric degenerations \cite{Manon,BGM-sympl} of a certain $(6g-6)$-dimensional manifold $\mathcal{N}_g$. These degenerations are associated to 3-valent graphs of genus $g$. Theorem 1 and the combinatorial non-abelian Torelli theorem \cite{BGM-torelli} imply that the monotone Lagrangian tori associated to different graphs are pairwise not Hamiltonian isotopic.

\section{Degenerations to singular symplectic spaces}
\subsection{Symplectic singular spaces}
K\"ahler varieties are complex analytic manifolds with a compatible symplectic structure. 
Any smooth complex projective variety $X$ or, more generally, complex submanifold of a K\"ahler variety $Y$ is a K\"ahler manifold itself.
Singular projective varieties are more delicate objects.
At the same time they are very useful as they might come (somewhat paradoxically) as degenerations of rigid (undeformable) K\"ahler manifolds such as $\cp^2$.

For the purposes of this paper it is convenient to keep the topological set-up as much as possible and the algebraic set-up as little as possible. 
With this in mind, we substitute singular projective varieties with the following rather general class of spaces generalizing complex algebraic varieties: {\em singular symplectic spaces in $\cp^N$}.
Furthermore, we may replace $\cp^N$ by another K\"ahler variety of sufficiently large dimension as follows.

Let $Y$ be a (smooth, but possibly non-compact) K\"ahler manifold of dimension $N$ with the symplectic form $\omega$.
Let $\Sigma\subset X\subset Y$ be closed subspaces such that 
\[
X^\circ=X\setminus\Sigma
\]
is a smooth connected complex submanifold of $Y$ of dimension $n$.
The form $\omega$ is a symplectic form on $X^\circ$. Its top exterior power $\omega^{n}$ is a volume form on $X^\circ$ and is defined everywhere on $Y$. Its cohomology class $[\omega^n]$ can be evaluated on any $2n$-cycle in $Y$.

\begin{defn}\label{def-sss}
We say that
$(X,\Sigma)$ is a {\em singular symplectic 
space} of (complex) dimension $n$ in the ambient K\"ahler manifold $Y$
if $\Sigma$ and $X$ are closed subspaces of $Y$  with $\Sigma\subset X$
such that the following conditions hold.
\begin{enumerate}
\item
The complement $X^\circ=X\setminus\Sigma$ is a smooth connected $n$-dimensional complex submanifold of $Y$.
\item
For any $x\in X^\circ$ the inclusion homomorphism $\phi_x:H_{2n}(X;\R)\to H_{2n}(X,X\setminus\{x\};\R)=H_{2n}(\C^n,\C^n\setminus\{0\};\R)=\R$ 
is an isomorphism. 

In particular, the fundamental class $[X]\in H_{2n}(Y;\R)$ is well-defined as the image of the fundamental class of $(\C^n,\C^n\setminus\{0\})$ under the composition of $\phi^{-1}$ and the homomorphism induced by the inclusion $X\subset Y$.
\item
We have
$$\int\limits_{X^\circ}\omega^n=[\omega^n][X],$$ where the right-hand side is defined by pairing of $[X]\in H_{2n}(Y;\R)$ and $[\omega^n]\in H^{2n}(Y;\R)$.

In particular, $X^\circ$ has a finite volume with respect to $\omega$.
\end{enumerate}
\end{defn}
We denote the restriction of the symplectic form from $Y$ to $X^\circ$ also with $\omega$.
It turns $X^\circ$ into a K\"ahler manifold. 

\begin{rmk}
Below by default we use \v Cech cohomology $H^k$ for our topological spaces in consideration. The corresponding homology groups (over fields $\Q$ and $\R$) are obtained as dual vector spaces of the corresponding cohomology groups. Note, however, that in most applications our main examples are CW-complexes, so in this case we may also equivalently choose singular or cellular cohomology as the principal homology theory.
\end{rmk}

\begin{defn}\label{def-qg}
A singular symplectic space $(X,\Sigma)$ is called a $\qg$-space if 
the following conditions hold.
\begin{enumerate}
\item The inclusion homomorphism 
\begin{equation}\label{eq-iota}
\iota: H^2(X;\R)\to H^2(X^\circ;\R)
\end{equation}
is injective.
\item The canonical class $K\in H^2(X^\circ;\R)$ is in the image of \eqref{eq-iota}.
\end{enumerate}
\end{defn}

Recall that $\bar a\in H^2(X;\R)$ is called a cohomological extension of $a\in H^2(X^\circ;\R)$ if $a$ is the image of $\bar a$ under
\eqref{eq-iota}.
The following proposition is simply a translation
of Definition \ref{def-qg} in the language of cohomological extension.
\begin{prop}
Suppose that $(X,\Sigma)$ is a $\qg$-space in $Y$.
\begin{enumerate}
\item
A cohomological extension of any class in $H^2(X^\circ;\R)$ is unique.
\item
A cohomological extension of the canonical class $K\in H^2(X^\circ;\R)$
exists.
\end{enumerate}
\end{prop}

For brevity, we denote the unique extension of the canonical class of a $\qg$-space $X$ also by $K\in H^2(X;\R)$.

\begin{rmk}
The second condition in Definition \ref{def-qg} mimics (in the context of symplectic singular spaces) the key property of $\Q$-Gorenstein algebraic varieties. \end{rmk}

The dense smooth part $X^\circ$
is naturally a symplectic $(2n)$-manifold with the symplectic structure $\omega$ induced by the
symplectic form
on $Y$.
The cohomological extension in $H^2(X;\R)$ of $[\omega]\in H^2(Y;\R)$ is given by the cohomology class of the symplectic form in $Y$.
Where confusion is not possible we denote it still with $[\omega]\in H^2(X;\R)$.

\begin{defn}\label{Xmonotone}
We say that a $\qg$-space $X\subset Y$ is {\em monotone} if the classes
$K$ and $[\omega]$ are proportional, i.e. $K=\alpha[\omega]\in H^2(X;\R)$ for some constant $\alpha\in \R$.
\end{defn}

\subsection{Lagrangian submanifolds of $\qg$-spaces}

Recall that for an orientable Lagrangian (or just only totally real) submanifold $L$ of a symplectic manifold $X^\circ$ we have a canonical trivialization along $L$ of the top exterior power $\Lambda^n(T^*X^\circ)$ of the cotangent bundle $T^*X^\circ$
(computed with respect to any symplectically compatible almost complex structure)
given by a real volume form on $L$. Even if $L$ is non-orientable, we get a trivialization of $\Lambda^n(T^*X^\circ)$ over any orientation preserving loop on $L$.

For a smooth compact surface $B\subset X^\circ$ with boundary such that  $\dd B\subset L$ is a collection of orientation-preserving loops in $L$, we define $K(B)\in\Z$ as the obstruction to extending this trivialization from $\dd B$ to $B$.
Since twice $\dd B$ is always orientation-preserving, for a general smooth closed surface $B\subset X^\circ$ with $\dd B\subset L$ we may consistently define
$K(B)\in\frac12\Z$ as one half of the obstruction to extending the corresponding trivialization on $2\dd B$. Clearly, the value $K(B)$ depends only on $[B]\in H_2(X^\circ,L)$.
The homomorphism
\begin{equation}
\label{mu-def}
\mu:H_2(X^\circ,L;\Q)\to\Q
\end{equation}
defined by $\mu([B])=-2K(B)$
is called the {\em Maslov index} of $[B]$.

\begin{prop}
If $(X,\Sigma)$ is a $\qg$-space and $L\subset X^\circ$ is a Lagrangian then the homomorphism \eqref{mu-def} factors through the inclusion homomorphism 
\begin{equation}\label{def-iota}
H_2(X^\circ,L;\Q)\to H_2(X,L;\Q).
\end{equation}
\end{prop}
\begin{proof}
The difference of two classes in $H_2(X^\circ,L;\Q)$ with the same image under \eqref{def-iota} comes from an absolute $\Q$-cycle in $X^\circ$ that is homologous in $X$ to a class supported by $L$.  By the second condition in Definition \ref{def-qg} the canonical class vanishes on this difference.
\end{proof}

\begin{coro}
If the homomorphism \eqref{def-iota} is surjective then 
there exists a unique homomorphism 
\begin{equation}
\label{qg-mu-def}
\mu:H_2(X,L;\Q)\to\Q
\end{equation}
extending \eqref{mu-def}.
\end{coro}
This homomorphism is still called the Maslov index.

Since $X\subset Y$ and $Y$ are K\"ahler, we may compute the symplectic area of a surface $B\subset X$ with boundary on a Lagrangian submanifold $L$. Since $L$ is Lagrangian,
the result $\omega(B)$ only depends on its homology class, thus producing the {\em area homomorphism}
\begin{equation}\label{def-omegaL}
\omega:H_2(X,L;\R)\to \R.
\end{equation}

\begin{defn}\label{XLmonotone}
Let $(X,\Sigma)$ be a $\qg$-space.
A Lagrangian submanifold $L\subset X^\circ$ is called {\em monotone in $X$} if \eqref{def-iota} is an epimorphism and the homomorphisms \eqref{def-omegaL} and \eqref{qg-mu-def} are proportional.
\end{defn}

Clearly, if $L\subset X$ is monotone then the $\qg$-space $X$ is also monotone, and the proportionality coefficients in Definitions \ref{Xmonotone} and \ref{XLmonotone} differ by multiplication by $-2$.

\subsection{Symplectic degenerations}
Let $0\in T\subset\C$ be the open unit complex disk (which we view as a symplectic surface of finite area),
$Y$ be a K\"ahler manifold (e.g. $\cp^N$), and $(\xx,\Sigma)$ be a singular symplectic space in $Y\times T$.
Denote by 
\begin{equation}\label{map-p}
p: \xx\to T
\end{equation}
the restriction of the projection $Y\times T\to T$ to $\xx$, and
by $X_t=p^{-1}(t)$ the fiber over $t\in T$. 

\begin{defn}
Let $(\xx,\Sigma)$ be a singular symplectic space in $Y\times T$ with $\Sigma\subset Y\times\{0\}$, such that
$X=\xx\cap (Y\times\{0\})$ is compact and $(X,\Sigma)$ is a singular space in $Y=Y\times\{0\}$.
We say that $X\subset Y$ is a {\em closed degeneration} of $X_t\subset Y=Y\times\{t\}$, $t\neq0$, if 
the following conditions hold.
\begin{enumerate}
\item
The map \eqref{map-p} is proper (in particular, each $X_t$ is compact).
\item
The restriction
$p|_{\xx^\circ}$ is a submersion.
\end{enumerate}
 In this case we say that  $X_t$, $t\in T\setminus\{0\}$, is a {\em symplectic smoothing} of $(X,\Sigma)$ and that $(\xx,p)$ is a {\em symplectic deformation} of $(X,\Sigma)$.
\end{defn}
We have $X^\circ=X\setminus\Sigma=\xx^\circ\cap X$.
The fiber $X=X_0$ is called the {\em central fiber}. 
Fibers $X_t$, $t\neq 0$, are called {\em general fibers}.

\begin{prop}\label{homeqX0}
If $(\xx,p)$ is a symplectic deformation then the inclusion of the central fiber $X\subset\xx$ induces an isomorphism in all cohomology groups $H^k(\xx;\R)\to H^k(X;\R)$
(and thus also in all homology groups $H_k(X;\R)\to H_k(\xx;\R)$).
\end{prop}
\begin{proof}
Recall that our basic homology theory is \v Cech cohomology and that $Y$ is a K\"ahler manifold. For homology computations, we may use coverings by small open disks in $Y$ and by products of small disks in $Y\times T$. Since $X$ is compact,
any such open covering of $X$ is induced by an open covering of $p^{-1}(U)$ with the same nerve, where $0\in U\subset T$ is a subdisk of $T$. But $p^{-1}(U)$ is a deformational retract of $\xx$ and thus the corresponding inclusion induces isomorphisms in homology.
\end{proof}

\begin{prop}\label{prop-sympl}
A smooth embedded path $\Gamma\subset T$ connecting $0$ and $t$ defines a closed set $\Sigma_\gamma\subset X_t$ of zero $\omega^n$-volume, and a symplectomorphism
\[
\Phi_\gamma:X^\circ\to X_t\setminus\Sigma_\gamma.
\]
\end{prop}
\begin{proof}
The complement $\tilde\Gamma=p^{-1}(\Gamma)\setminus\Sigma$ is a smooth $(2n+1)$-submanifold of $Y\times T$ with boundary composed of $X^\circ$ and $X_t$. The form $\omega$ is symplectic on each $X_s$, $s\in T$. The annihilator of the restriction of $\omega$ to the interior of $\tilde\Gamma$ gives the characteristic foliation ${\mathcal F}$ defining a connection on 
\[
p|_{\tilde\Gamma}:\tilde\Gamma\to\Gamma
\]
that preserves the symplectic form.
The leaf of ${\mathcal F}$ passing through a point $x_0\in X^\circ$ is a  closed interval whose other end is a point $x_t\in X_t$.
We define $\Phi_\gamma(x_0)=x_t$ and $\Sigma_\gamma=X_t\setminus\Phi_\gamma(X^\circ)$.
The locus $\Sigma_\gamma$ is clearly closed. It remains to prove that $\int\limits_{X_t\setminus\Sigma_\gamma}\omega^n=\int\limits_{X_t}\omega^n$.
But the latter integral coincides with the evaluation of cohomological class $[\omega^n]$ (which can be represented in \v Cech cohomology with the help of De Rham isomorphism) at $[X_t]$. By Proposition \ref{homeqX0} $[X_t]\in H_{2n}(\xx;\R)$ is in the image of the inclusion homomorphism from $H_{2n}(X;\R)$. With the help of a ($\xx^\circ$-valued) section of $p$ we see that $[X_t]$ has to be the image of $[X]$.
Conservation of symplectic form by $\Phi_\gamma$ implies that $$\int\limits_{X_t\setminus\Sigma_\gamma}\omega^n=\int\limits_{X\setminus\Sigma}\omega^n=[\omega^n][X]=[\omega^n][X_t].$$
\end{proof}

The following statement can be formally obtained from Proposition \ref{prop-sympl} by restricting to a subdisk $s,t\in T'\subset T\setminus\{0\}$ and translating it by $-s$, so that the smooth fiber $X_s$ becomes the new central fiber with $X^\circ_s=X_s$.
\begin{coro}\label{coro-sy}
Any two general fibers $X_t, X_s$, $s,t\in T$ are symplectomorphic.
\end{coro}

\begin{defn}
A symplectic deformation $(\xx,p)$ in $Y\times T$ is called a $\qg$-deformation if both $(\xx,\Sigma)$ and $(X,\Sigma)$ are $\qg$-spaces.
Accordingly, we say that $X_t$ is a {\em $\qg$-smoothing} of $(X,\Sigma)$.
\end{defn}
The main algebro-geometric source of $\qg$-deformations are
$\Q$-Gorenstein deformations in the sense of \cite{KSB}.

Suppose that $A\subset X^\circ$ is a subset.
\begin{defn}
We say that $(X_t,A_t)$ is a $\qg$-smoothing of $(X,A)$ if $X_t$ is a $\qg$-smoothing of $(X,\Sigma)$, and there exists a path $\Gamma\subset T$ as in Proposition \ref{prop-sympl} such that $A_t=\Phi_\gamma(A)$.
\end{defn}

\begin{prop}
Suppose that $(\xx,p)$ is a symplectic deformation of a $\qg$-space $(X,\Sigma)$.
Suppose that $H_1(X_t;\R)=0$ for 
a general fiber $X_t$, $t\in T\setminus\{0\}$.
Then $(\xx,p)$ is a $\qg$-deformation if and only if the canonical class $K_{X_t}\in H^2(X_t;\Q)$ vanishes on the kernel of the inclusion homomorphism
\begin{equation}\label{Xthom}
H_2(X_t;\R)\to H_2(\xx;\R)\approx H_2(X;\R).
\end{equation}
\end{prop}
\begin{proof}
If the value of $K_{X_t}$ is not zero on an element of the kernel of \eqref{Xthom} then the value of $K_{\xx^\circ}$ is also not zero on this element, and thus $K_{\xx^\circ}$ is not extendable to $\xx$.

Consider the converse direction, i.e. suppose that $K_{X_t}$ vanishes on the kernel of \eqref{Xthom}.
Note that since $\xx^\circ\supset X^\circ$ and $(X,\Sigma)$ is a $\qg$-space, the homomorphism 
\begin{equation}\label{surjXX}
H_2(\xx^\circ;\R)\to H_2(\xx;\R)\approx H_2(X;\R)
\end{equation}
is surjective. The isomorphism here is provided by Proposition \ref{homeqX0}. Thus we are left to prove that $K_{\xx^\circ}$ vanishes on the kernel of \eqref{surjXX}. Let $S\subset\xx^\circ$ be a closed smooth surface representing an element $\alpha$ of this kernel. We may assume that it is transversal to $X^\circ\subset\xx^\circ$.
Since $p(S)\subset T$ is compact, the intersection number of $X^\circ$ and $S$ in $\xx^\circ$ vanishes, and we may assume (after adding some handles to $S$) that $S\cap X=\emptyset$.

Consider a radius $R\subset T$, i.e. a half-open smooth path connecting $0$ and the boundary of $T$. We may assume that $S$ and $p^{-1}(R)$ are transverse, and thus intersect along a smooth curve in $p^{-1}(R)\approx X_t\times R$. Since $H_1(X_t;\R)=0$ we may find a representative of $\alpha=[S]$ contained in $\xx^\circ\setminus p^{-1}(R)\approx X_t\times (T\setminus R)$ and thus in $X_t$. The value of $K_{X_t}$ and thus of $K_{\xx^\circ}$ on this representative vanishes by assumption.
\end{proof}

\begin{coro}\label{qg-auto}
Suppose that $(\xx,p)$ is a symplectic deformation of a $\qg$-space $(X,\Sigma)$ with general fiber $X_t$.
If $H_1(X_t;\R)=0$ and $H_2(X_t;\R)\approx H_2(X;\R)$ then $(\xx,p)$ is a $\qg$-deformation.
\end{coro}

\begin{prop}\label{Lmonotone}
Suppose that $(X_t,L_t)$ is a $\qg$-smoothing of $(X,L)$ where $L\subset X^\circ$ is a Lagrangian submanifold monotone in $X$.
Then $L_t\subset X_t$ is a monotone Lagrangian submanifold. 
\end{prop}
\begin{proof}
Given a membrane $B_t\subset X_t$, $\dd B_t\subset L_t$,
we get a new membrane $B\subset\xx^\circ$, $\dd B\subset L$,
by attaching to $B_t$ a cylinder $\Gamma\times\dd B_t$ along the symplectic connection defined by the path $\Gamma$. By construction we have $\omega(B)=\omega(B_t)$.
The trivializations of the canonical classes of $X^\circ$ and $X_t$ over (the orienting double cover of) $L$ and $L_t$ define the trivialization of the canonical class of $\xx^\circ$ over  (the orienting double cover of) $L\cup L_t$. Thus we may speak of the values of the relative canonical class
and we get 
\begin{equation}\label{muB0}
\mu_{\xx^\circ}(B_t)=\mu_{\xx^\circ}(B).
\end{equation}
But these values depend only on the classes of $B_t$ and $B$ in $H_2(\xx,L_t;\Q)$ and $H_2(\xx,L;\Q)$ since $\xx$ is a $\qg$-space.
Furthermore, since $B_t\subset X_t$ and $X_t$ is a regular fiber of $p$, we have $\mu_{\xx^\circ}(B_t)=\mu_{X_t}(B_t)$.

To conclude the proof we note that the membrane $B$ must be homologous in $H_2(\xx,L)$ to a membrane $B_0$ contained in $X$, since the embedding $X\subset \xx$ is a homotopy equivalence. We get $\mu_{\xx^\circ}(B_0)=\mu_{X}(B_0)$ since $X$ is a $\qg$-space.
Thus $\mu_{X}(B_0)=\mu_{X_t}(B_t)$. Since $\omega(B_0)=\omega(B)=\omega(B_t)$ we conclude that $L_t\subset X_t$ is monotone with the same monotonicity coefficient as $L\subset X$.
\end{proof}

\begin{rmk}
For monotonicity of $L_t$ the hypothesis that $(\xx,\Sigma)$ is a $\qg$-space is crucial.
While $\Phi_\gamma(L)\subset X_t$ is Lagrangian for any symplectic smoothing $X_t$
of a symplectic singular space $(X,\Sigma)$, it does not have to be monotone even if  $(X,\Sigma)$ is a $\qg$-space and $L\subset X^\circ$ is a monotone unless we assume that $X_t$ is a $\qg$-smoothing, i.e. $(\xx,\Sigma)$ is also a $\qg$-space. \end{rmk}

\section{Toric degenerations}
In this section we apply the framework described in the previous section to the case when the ambient space $Y=\cp^N$ is projective, $X=X\subset Y$ is a normal toric variety whose embedding to $Y$ is defined by a convex lattice polyhedron $\Delta\subset\R^n$ and $\Sigma\subset X$ is the singular locus of $X$. In this case smoothings of $X$ can be produced in the framework of the patchworking construction introduced by Viro \cite{Vi-patch} (see also its interpretation in \cite{GKZ} and \cite{ItMiSh}). We plan to give a detailed description of applying patchworking in the context of symplectic smoothings in our future paper \cite{GaMi}.

\subsection{Toric geometry set-up} Here we recall the basic language of toric geometry (cf. e.g. \cite{CdS} and the references therein). Let $M\approx\Z^n$ be an integer lattice of rank $n$ and $N=M^*=\operatorname{Hom}(M,\Z)$ its dual lattice. Denote by $\mr=M\otimes\R$ the corresponding real vector space and by $\nr=N\otimes\R=\mr^*$ its dual space. Let $$\Delta\subset\nr$$ be a convex polyhedron with integer vertices and non-empty interior. Denote by $\ff\subset M\setminus\{0\}$ the set of the primitive outwards supporting vectors for the facets of $\dd\Delta$.

The polyhedron $\Delta$ determines a toric K\"ahler manifold $X=X_\Delta$. There are several constructions useful in different frameworks. As a symplectic manifold, $X$ may be obtained from the quotient $\Delta\times (\mr/M)$ of the cotangent bundle $T^*\Delta=\Delta\times\mr$ by the lattice $M\subset\mr$. Namely, over a face $\Delta'\subset\dd\Delta$ 
we take a quotient of the fiber torus $\mr/M$ by the subtorus annihilating the vector subspace of $\nr$ parallel to $\Delta'$.
The projection to $\Delta$ yields the {\em moment map}
\[
\mu_\Delta:X_\Delta\to\Delta.
\]
For $x\in\Delta$ the inverse image $\mu_\Delta^{-1}(x)$ is a real torus of dimension equal to the dimension of the face of $\Delta$ containing $x$.

As a complex manifold $X$ can be obtained with the help of $\ff$. Namely, each face $\Delta'\subset\dd\Delta$ gives us a subset $\ff_{\Delta'}$ consisting of the supporting vectors to the facets whose closure contains $\Delta'$.
The elements of $\ff_{\Delta'}$ generate a cone in $M$, with the dual cone $\sigma_{\Delta'}\subset N$.
The spectrum of the semigroup algebra (over  $\C$) of $\sigma_{\Delta'}$ is a (possibly singular) complex chart $U_{\Delta'}$. The complex variety $X$ is obtained by gluing these charts. All the charts contain the complex torus $\mc\approx(\C^\times)^n$.
If $\ff_{\Delta'}$ is an integer basis of the subspace $\mr^{\Delta'}\subset\mr$ it spans,
then the chart $U_{\Delta'}$ is isomorphic to
$$\tilde U_{\Delta'}=(\C^\times)^{\dim\Delta'}\times\C^{n-\dim\Delta'}.$$ 
If it is only a rational basis then the chart $U_{\Delta'}$ is isomorphic to the quotient of $\tilde U_{\Delta'}$ by a finite group (with the number of elements equal to the discriminant of $\ff_{\Delta'}$). Considered together, symplectic and complex constructions yield the following proposition.
\begin{prop}\label{toric-way2}
The singular locus of $X_\Delta$ is given by 
\[
\Sigma=\bigcup\mu_\Delta^{-1}(\Delta'),
\]
where the union is taken over faces $\Delta'\subset\dd\Delta$ such that $\ff_{\Delta'}\subset M$ is not an integer basis of $\mr^{\Delta'}$.

If $\Delta$ is simple, i.e. each of its vertices is adjacent to $n$ facets, then a neighborhood of each singular point of $\Sigma$ in $X$ can be presented as the quotient singularity of $\C^n$ by a finite group.
\end{prop}
In particular, $X_\Delta$ is smooth outside of a complex codimension 2 locus and $H_{2n}(X;\Q)=\Q$.

There is yet the third toric variety construction exhibiting $X_\Delta$ as a subvariety of $\cp^m$ with $m=\#(\Delta\cap N)-1$, thus giving its non-singular part a K\"ahler structure. For this purpose consider a bijection between the $m+1$ homogeneous coordinates of $\cp^m$ and the lattice points of $\Delta$.  Each lattice point of $\Delta$ is a multiplicatively linear functional on $M_{\C}=M\otimes\C^\times\approx (\C^\times)^n$. Taken together, they define a map to $(\C^\times)^{m+1}$ and, after projectivization, to $\cp^m$. The closure $X_\Delta\subset\cp^m$ of the image recovers at the same time the symplectic toric variety of the first construction and the complex toric variety from the second construction,  
see \cite{GKZ}. This implies the following proposition.
\begin{prop}
For each lattice polytope $\Delta\subset\nr$ the pair $(X_\Delta,\Sigma)$ is a singular symplectic space.
\end{prop}
The smooth locus $X^\circ=X\setminus\Sigma$ is a K\"ahler manifold.
The inverse image \[ Z=\mu^{-1}_\Delta(\dd\Delta)\] is a complex (reducible) hypersurface. Smooth points of $Z$ are simple poles of the meromorphic volume form 
\begin{equation}\label{formOmega}
\Omega=\frac1{(2\pi i)^n}\bigwedge\limits_{j=1}^n \frac{dz_j}{z_j}.
\end{equation}
Thus the homology class (with closed support) of $Z\cap X^\circ$ is dual to minus the canonical class on $X^\circ$. Here $z_j$, $j=1,\dots,n$, are basis coordinates in the complex torus $U_\Delta=\mc\approx(\C^\times)^n$.

Each element $v\in M$ can be viewed by duality as a functional $v:\nr\to\R$.
For a facet $\delta\subset\Delta$ denote with $v_\delta\in\ff\subset M$ its outward primitive integer supporting vector. Then $v_\delta|_\delta$ is constant. Denote $h_\delta=v_\delta(\delta)$.
The condition that $0$ is contained in the interior of $\Delta$ is equivalent to the condition that $h_\delta>0$ for all facets $\delta\subset\Delta$.
In this case the fiber
\[
L=\mu^{-1}_{\Delta}(0)\subset X
\]
is a Lagrangian torus.
Furthermore, the restriction of $\Omega$ to $L$ is a real volume form (of total volume 1). Thus
the Maslov index $m(B)$ of a membrane $B\subset X$, $\dd B\subset L$, coincides with twice the intersection number $\#(B\cap Z)$,
\begin{equation}\label{mu-index}
m(B)=2\#(B\cap Z).
\end{equation}

\begin{lem}\label{lem-mono}
The Lagrangian submanifold $L\subset X^\circ$ is monotone if and only if $h=h_\delta$ does not depend on the choice of the facet $\delta\subset\Delta$. 
\end{lem}
\begin{proof}
Connect $0$ and a point inside $\delta$ with a simple path $\gamma\subset\Delta$. 
Let $D_\gamma\subset X$ be the disk obtained from $\gamma\times S\subset\Delta\times(\mr/M)$
after the quotient. Here $S\subset \mr/M$ is a simple closed curve parallel to $v_\delta$.
By \eqref{mu-index} we have $m(D_\gamma)=2$. But $\omega(D_\gamma)=h_\delta$ by the construction of the symplectic form $\omega$ on the cotangent space $T^*\Delta=\Delta\times(\mr/M)$. The lemma follows from the observation that the disks $D_\gamma$ generate the group $H_2(X^\circ,L;\Q)$.
\end{proof}

\subsection{Toric $\qg$-spaces and Maslov indices.}
In general, $(X_\Delta,\Sigma)$ does not have to be a $\qg$-space. 
\begin{prop}
A pair $(X_\Delta,\Sigma)$ is a $\qg$-space if and only if for each vertex $\Delta'\subset\Delta$ (a face of dimension 0) the set $\ff_{\Delta'}\subset M$ is contained in an affine hyperplane.
\end{prop}
\begin{proof}
To compute the group $H_2(X_\Delta;\R)$ in terms of $\Delta$ we 
consider the vector space $V_\Delta$ whose basis is in 1-1 correspondence with the facets of $\Delta$ and the linear map $\lambda:V_\Delta\to\mr$ sending each facet to its outward normal vector. For the dense torus $U_\Delta\approx(\C^\times)^n$ in $X_\Delta$ we have
$$H_2(X_\Delta,U_\Delta;\R)=V_\Delta.$$ The long exact sequence of $(X_\Delta,U_\Delta;\R)$ implies that $H_2(X_\Delta;\R)=\ker\lambda$.
Let $\operatorname{Sk}^{n-2}(\Delta)$ be the union of the faces of $\Delta$ of dimension up to $n-2$.
Since the inclusion homomorphism $$H_2(X\setminus\mu_\Delta^{-1}(\operatorname{Sk}^{n-2}(\Delta));\R)\to H_2(X;\R)$$ is surjective, the inclusion  homomorphism $H_2(X^\circ;\R)\to H_2(X;\R)$ is surjective as well.

Since $\Delta'$ is a vertex, its inverse image $\mu^{-1}(\Delta')$ is a point. The chart  $U_{\Delta'}$ is contractible and contains a fiber torus $L_x=\mu^{-1}(x)$, $x\in \Int\Delta$. 
Connecting $x$ to a facet $\delta$ adjacent to $\Delta'$ with a simple path we get a disk $D_\gamma$ of Maslov index 2 as in
the proof of Lemma \ref{lem-mono}. The class $$
[D_\gamma]\in H_2(X_\Delta,L_x;\R)=H_2(X_\Delta,U_\Delta;\R)=V_\Delta
$$
is the basis vector corresponding to $\delta$.
Thus $K_{X^\circ}$ has equal values on all basis vectors in $V_\Delta$.
But $D_\gamma$ is contained in $U_{\Delta'}$. 
Contractibility of $U_{\Delta'}$ implies
$$H_2(U_{\Delta'},L_x)=H_1(L_x)=M\approx\Z^n.$$
If $K_{X^\circ}$ is extendable to $X^\circ\cup U_{\Delta'}$ then any linear relation on the elements of $\ff_{\Delta'}$ has to have zero sum of coefficients. Thus the elements of  $\ff_{\Delta'}$ must be contained in a hyperplane of $\mr$. Conversely, if they are contained in such a hyperplane, we may find an extension of $K_{X^\circ}$ to $X^\circ\cup U_{\Delta'}$, and thus also to $X$.
\end{proof}
The hyperplane containment condition is vacuous in the case when each vertex of $\Delta$ is adjacent to $n$ facets.
\begin{coro}\label{coro-sim}
If $\Delta$ is simple then $(X_\Delta,\Sigma)$ is a $\qg$-space.
\end{coro}

From now on we suppose that $0\in\Delta\subset\nr$ is a simple monotone convex lattice polyhedron, so that $L=L_0$ is a monotone Lagrangian torus in the $\qg$-space $(X_\Delta,\Sigma)$.

Suppose that $B\subset X_\Delta$ is a membrane, i.e. the image of a smooth surface map in $X$ such that $\dd B\subset L$.
Suppose that $x\in B\cap Z$ is isolated.
By Proposition \ref{toric-way2}, a neighbourhood of $x$ is locally isomorphic to the quotient of $\C^n$ by a group of order $l\in\Z_{>0}$.
We define the local intersection number $\#_x(B\cap Z)$ as $\frac 1l$ times the local intersection number of the inverse images of $B$ and $Z$ in $\C^n$. It is a rational number.

Recall that $U_\Delta=\mu^{-1}_\Delta(\Int\Delta)$ has the canonical complex structure of $\mc\approx(\C^\times)^n$. Thus it comes with the natural map
\[
\Log:(\C^\times)^n\approx U_\Delta\to\mr\approx \R^n,
\]
given in coordinates by $\Log(z_1,\dots,z_n)=(\log|z_1|,\dots,\log|z_n|)$.

A non-zero element $v\in M$ generates a subgroup in $\mc$, and thus also a holomorphic cylinder $C_v\subset (\C^\times)^n$ such that $\Log(C_v)\subset\mr$ is the line passing through the origin and parallel to $v$.
We refer to any multiplicative translate $a C_v$, $a\in\mc$ as a {\em holomorphic cylinder parallel to $v$}.

By a {\em proper holomorphic curve in $X$ with boundary in $L$} we mean the image $B=f(S)\subset X$ of a non-constant continuous map $f:S\to X$ from a compact Riemann surface $S$ with boundary $\dd S$ such that $f|_{S\setminus\dd S}$ is holomorphic and $f(\dd S)\subset L$. We say that $B$ is irreducible if $S$ is connected. We refer to an irreducible proper holomorphic curve with boundary in $L$ (as well as to its image in $X$) as a {\em holomorphic membrane} for $L$ if its {\em boundary} $\dd B=f(\dd S)$ is non-empty.
\begin{prop}
Suppose that $B\subset X_\Delta$ is 
a holomorphic membrane for $L_x$, $x\in\Int\Delta$.
Then $B\cap Z$ is a finite set and 
the Maslov index of $B$ can be computed as
\begin{equation}\label{eq-min}
m(B)=2\sum\limits_{z\in B\cap Z}\#_z(B\cap Z),
\end{equation}
and $\#_z(B\cap Z)>0$ for any $z\in B\cap Z$.
\end{prop}
\begin{proof}
Since $Z$ is disjoint from a neighborhood of $\dd B\subset L$, the set of isolated points of $Z\cap B$ must be finite, $l(\#_z(B\cap Z))$ is the local intersection number of a curve and a hypersurface in $\C^n$ and thus is a positive integer number.
The formula \eqref{eq-min} is a corollary of \eqref{mu-index} and the surjectivity of the homomorphism $H_2(X^\circ,L)\to H_2(X,L)$.
\end{proof}

\begin{lem}\label{lem-const}
Suppose that  $B=f(S)\subset X_\Delta$ is
a holomorphic membrane for $L_x=\mu_\Delta^{-1}(x)$, $x\in\Int\Delta$,
such that $f^{-1}(Z)\subset S$ consists of a single point.

Then $\nu_{\C}|_B$ is constant for any multiplicatively linear map 
$$\nu_{\C}:\mc\to\C^\times$$ defined as the pairing with an element
$\nu\in N$ such that $\nu([\dd B])=0$.
\end{lem}
Here we view $[\dd B]$ as the element of $H_1(L_x;\R)=\mr$ (sitting on the lattice $M\subset\mr$).
\begin{proof}
Note that $\nu_{\C}:\mc\to \C^\times$ is a product of integer powers of the holomorphic coordinates in the dense torus $U_\Delta\subset X_\Delta$.
Let $w\in S\setminus \dd S$ be the (only) point of $f^{-1}(Z)$. Its image $f(w)\in Z$ is contained in the toric chart $U_{\Delta'}$ for a face $\Delta'\subset\Delta$, and thus the meromorphic function $\nu_{\C}\circ f$ near $w$ can be expressed in a local coordinate $z$ near $w$ as
$$c(z-w)^k+o((z-w)^k),$$
$c\in\C^\times$, $k\in\Z$.
If $k>0$ then a small punctured disk centered at $w$ covers a small disk in $\C$ centered at 0. If $k<0$, we get such covering for $1/(\nu_{\C}\circ f)$. In both cases we conclude that the linking number of $\nu_\C\circ f(\dd S)$ and $0$ in $\C$ is $k\neq 0$. 
We get a contradiction since this number is given by $\nu[\dd B]=0$.
Thus $k=0$ and
the function $\nu_{\C}\circ f$ can be extended to a continuous function on $S$ holomorphic on $S\setminus\dd S$ and constant on $\dd S$ (since $f(\dd S)\subset L_x$).
By the maximum principle, $\nu_{\C}\circ f$ is constant.
\end{proof}

\begin{coro}\label{coro-point}
Suppose that  $B=f(S)\subset X_\Delta$ is a holomorphic membrane for
$L_x=\mu_\Delta^{-1}(x)$, $x\in\Int\Delta$, such that $f^{-1}(Z)\subset S$ consists of a single point.

Then $[\dd B]\neq 0\in H_1(L)=M$, and $B\cap (\mc)$ is contained in a holomorphic cylinder $aC_{[\dd B]}$, $a\in (\C^\times)^n$,  parallel to $[\dd B]\in M$.
\end{coro}
\begin{proof}
By Lemma \ref{lem-const}, $\nu_{\C}$ is constant on the set $B\cap (M\otimes\C^\times)$ for any $\nu\in N$ with $\nu[\dd B]=0$.
Thus if $[\dd B]=0$ then $B$ must consist of a single point which leads us to a contradiction.
If $[\dd B]\neq 0$ then the multiplicative projection of $B\cap (M\otimes\C^\times)$ along $[\dd B]$ consists of a single point in $(\C^\times)^{n-1}$, and therefore, $B\cap (M\otimes\C^\times)$ is contained in $aC_{[\dd B]}$ for some $a\in(\C^\times)^n$.
\end{proof}

\section{Holomorphic membranes in toric $\qg$-degenerations}
As it was suggested in \cite{EP}, monotone Lagrangian submanifolds may be symplectically distinguished by the enumeration of holomorphic disks of the smallest possible Maslov index. The disks are holomorphic with respect to an almost complex structure compatible with the symplectic structure.
In the case of orientable monotone Lagrangian with the positive proportionality coefficient between the Maslov index and the area, any holomorphic membrane has the Maslov index of at least 2, so no bubbling of disks is possible.
By the Riemann-Roch formula, the expected dimension of the space of Maslov index 2 disks is $n-1$.
Thus by the Gromov compactness \cite{G85} we have finitely many of such disks passing through a fixed point of the Lagrangian in the case if these disks are regular.

Consider a closed symplectic degeneration $p:\xx\to T$, $X_t=p^{-1}(t)$, $t\in T\ni 0$.
Let $L=L_0\subset X_0$ be an orientable Lagrangian submanifold (disjoint from the singular locus $\Sigma\subset X_0\subset\xx$).
The image $L_t=\Phi_\gamma(L_0)$ of $L_0$ under the symplectomorphism $\Phi_\gamma:X^\circ\to X_t\setminus\Sigma_\gamma$ provided by Proposition \ref{prop-sympl} is an orientable Lagrangian in the smooth symplectic manifold $X_t$. If $\xx$ is a $\qg$-degeneration and $L_0\subset X_0$ is monotone then by Proposition \ref{Lmonotone} the Lagrangian $L_t\subset X_t$ is also monotone.
In other words, $\qg$-smoothings of $(X,L)$ where $L$ is a monotone Lagrangian submanifold of $X$ produce monotone Lagrangian submanifolds $L_t$ diffeomorphic to $L=L_0$ in the smooth symplectic manifold $X_t$, $t\neq 0$.

For the rest of the paper we assume that $(X_t,L_t)$ is a $\qg$-smoothing of $(X_0,L_0)$ where $L_0\subset X_0$ is a monotone orientable Lagrangian.
Thus for any relative homology class $\beta_t\in H_2(X_t,L_t)$ with $m(\beta_t)=2$ the number of holomorphic disks $\nd(\beta_t)\in\Z_2$ passing through a point $a\in L_t$ for a generic compatible almost complex structure is well-defined.
It is invariant under symplectomorphisms, and does not depend on the choice of a generic compatible almost complex structure.
Here we reduce the number of disks mod 2 to avoid a need to discuss orientation issues.

Since $L\subset X^\circ$, we have a canonical identification $H_1(L_t)=H_1(L)$ for any $t\in T$ by  Proposition \ref{prop-sympl}.
For a class $\alpha\in H_1(L;\R)$ we define
\begin{equation}\label{alpha-disks}
\nd(\alpha)=\sum\limits_{\beta_t} \nd(\beta_t)\in\Z_2,
\end{equation}
where the sum is taken over all $\beta_t\in H_2(X_t,L_t)$ of Maslov index 2 and such that $\dd\beta_t=\alpha
\in H_1(L_t;\R)$.
By the Gromov compactness and monotonicity, the sum on the right-hand side is finite since the area of all holomorphic disks in the sum is the same (and thus bounded).
By Corollary \ref{coro-sy}, the sum does not depend on the choice of $t\neq 0$.
We define
\begin{equation}
\ddd(L)=\operatorname{ConvexHull}\{\alpha\in H_1(L;\R)\ |\ \nd(\alpha)\neq 0\}.
\end{equation}
This is a convex lattice polygon in the vector space $H_1(L;\R)$ with the lattice $H_1(L)/\operatorname{Tor} H_1(L)$.
Any symplectomorphism of $(X_t,L_t)$ induces a lattice-preserving isomorphism of $H_1(L;\R)$. As a lattice polygon, $\ddd(L)$ is a symplectic invariant of $(X_t,L_t)$.

Following \cite{NNU}, \cite{NNU2}, we consider the case when $L=\mu_\Delta^{-1}(0)$ is a monotone Lagrangian torus in the toric variety $X_\Delta$ which is a symplectic $\qg$-space. Note that by Corollary \ref{coro-sim} the latter holds for any simple lattice polygon $\Delta\subset\R^n$. More generally, the $\qg$-condition coincides with the $\Q$-Gorenstein condition on the normal toric variety $X_\Delta$. Denote
$$\Delta^*=\operatorname{ConvexHull}(\ff)\subset\mr,$$
where $\ff\in\mr$ is the set of primitive outward supporting vectors for the facets of $\Delta$.

\begin{thm}\label{thm-main}
For any $\qg$-smoothing $(X_t,L_t)$ of $(X_\Delta,L)$ we have
$$\ddd(L)=\Delta^*.$$
\end{thm}
\begin{proof}
Suppose that $\alpha\in H_1(L;\R)$ is a class such that $\nd(\alpha)\neq 0$.
Consider the smooth embedded path $\Gamma\subset T$ such that $L_t=\Phi_\gamma(L)$. Let $\Gamma'\subset T$ be a smooth continuation of $\Gamma$ beyond the endpoint, so that $\Gamma'$ is diffeomorphic to an open interval.
Since each point of $\Gamma'$ is connected to $0$ by a subpath of $\Gamma'$, the space $X_s$ and its Lagrangian submanifold is well-defined for all $s\in\Gamma'$.
For each $s\in\Gamma\setminus\{0\}$ there exists a holomorphic disk membrane in $(X_s,L_s)$ whose boundary represents $\alpha$.
All these disks are contained in the symplectic product $Y\times S^2$, $S^2\supset T$, and have the same symplectic area.
The boundaries of the disks are contained in 
$\bigcup_{s\in\Gamma'}L_s$. This is a non-proper isotropic submanifold of $Y\times T$ which is contained in a non-proper Lagrangian submanifold $\mathcal L$ of $Y\times T$.
Furthermore, the boundaries of the disks are contained in a compact subset of  $\mathcal L$.
By the Gromov compactness we may extract a subsequence of holomorphic disk membranes converging to a possibly reducible (bubbled) curve $B\subset X$ of Maslov index 2 whose boundary represents $\alpha$.

The components of $B$ consist of holomorphic disk membranes $B_j=f_j(S_j)$ for $(X_\Delta,L)$ as well as closed rational curves $A_k\subset X_\Delta$. We define the Maslov numbers of the closed components $A_k$ by $m(A_k)=-2K.[A_k]$. The Maslov numbers of $A_k$ and $B_j$ are positive rational numbers adding to $m(B)=2$.
Suppose that $C=f_C(S_C)$ is a component of $B$ with the property that $C\setminus Z\neq 0$ (recall that $Z$ is union of toric divisors given by \eqref{formOmega}). Any point $z\in f^{-1}_C(Z)$ is isolated. Furthermore, the sum of the Maslov numbers of all components $A_k$ bubbling at $z$ and twice the local intersection number of $f_C$ and $Z$ at $z$ is a positive even number. Indeed the half of this number coincides with the obstruction to extending the trivialization of $\Lambda^n(T^*(X_\Delta\setminus\Sigma))$ given by \eqref{formOmega} from a small positive loop around $z$ to a disk around $z$.
Thus $C\subset X_\Delta$ must be a disk component with $[\dd C]=\alpha$ and $f_C^{-1}(Z)$ consists of a single point $z$.

The image $\mu_\Delta(z)$ is contained in the relative interior of a face of $\Delta$.
If $\Delta'$ is any vertex of this face then $C\subset U_{\Delta'}$.
Since the vertex chart $U_{\Delta'}$ is contractible, we have 
$$H_2(U_{\Delta'},L;\R)= H_1(L;\R)=\mr.$$
The Maslov index of membranes in $(U_{\Delta'},L)$ yields a linear function $\mr\to\R$ taking value 2 at all elements of $\ff_{\Delta'}$ (spanning a hyperplane in $\mr$ by the $\qg$-condition).
Since $z\in U_{\Delta'}$, the vector $\alpha\in\mr$ belongs to the cone spanned by $\ff_{\Delta'}$.
If $m(C)=2$ then $\alpha$ belongs to the convex hull of the set $\ff_{\Delta'}$.
If $0<m(C)<2$ then $\alpha$ belongs to the convex hull of the union of $\ff_{\Delta'}$ and the origin.
This implies $\ddd(L)\subset\Delta^*$.

For the opposite inclusion it suffices to prove that $\nd(\alpha)\neq 0$ if $\alpha$ is a vertex of $\Delta^*$.
Consider the holomorphic cylinder $aC_\alpha\subset M\otimes\C^\times$ passing through $a\in L$.
The complement $aC_\alpha\setminus L$ consists of two components.
Let $B$ be the component whose boundary orientation at $L$ agrees with $\alpha$.
The closure $\bar B\supset B$ in $X_\Delta$ is a regular holomorphic disk intersecting $Z$ transversally at a single point $w$ (in the divisor corresponding to $\alpha$).
We have $m(\bar B)=2$ by \eqref{mu-index}.
Furthermore, for any other (bubbled) holomorphic limit $C\cup\bigcup\limits_kA_k$ of holomorphic disks in $(X_s,L_s)$ its membrane part $C$ must intersect $Z$ transversally at the same point $w$ by Corollary \ref{coro-point}.
Thus $m(C)=2$ and $C=\bar B$, so no bubbling can appear and $\nd(\alpha)=1$.
\end{proof}

Let $\Delta,\Delta'\subset \nr\approx\R^n$ be two monotone lattice polytopes and $\Delta^*,(\Delta')^*\subset\mr\approx\R^n$ be their dual polytopes.
Let $X_{\Delta}$, $X_{\Delta'}$ be the corresponding toric varieties with the singular loci $\Sigma$, $\Sigma'$ and the monotone Lagrangian tori $L=\mu_\Delta^{-1}(0)\subset X_\Delta$, $L'=\mu_{\Delta'}^{-1}(0)\subset X_{\Delta'}$.
Suppose that $(X_\Delta,\Sigma)$ and $(X_{\Delta'},\Sigma')$ admit $\qg$-smoothings $X_t$ and $X'_t$.
Let $L_t\subset X_t$ and $L'_t\subset X'_t$ be the monotone Lagrangian tori given by Proposition \ref{Lmonotone}.

\begin{corollary}
If the pairs $(X_t,L_t)$ and $(X'_t,L'_t)$ are symplectomorphic then the polytopes $\Delta^*,(\Delta')^*\subset\mr$ coincide up to an automorphism of the lattice $M$.
\end{corollary}
\begin{proof}
The automorphism of $M=H_1(L_t)=H_1(L'_t)$ is induced by the restriction to $L_t$ of the symplectomorphism between $X_t$ and $X'_t$.
\end{proof}

In other words, if $X_\Delta$ and $X_{\Delta'}$ are different toric degenerations of the same symplectic manifold then the tori $L_t$ and $L'_t$ are symplectically different.
\begin{corollary}
The number of non-symplectomorphic Lagrangian tori in a complex projective variety is bounded from below by the number of its $\Q$-Gorenstein toric degenerations.  
\end{corollary}

\begin{exa}\label{exa-dim2}
In the case $n=2$ there is an infinite number of different $\Q$-Gorenstein toric degenerations of Del Pezzo surfaces, see \cite{HaPr} for $\dim H_2(X_\Delta;\R)=1$ and in \cite{CMG} for $\dim H_2(X_\Delta;\R)>1$.
This implies that the number of distinct monotone Lagrangian tori in $\cp^2$ and all other Del Pezzo surfaces is also infinite (cf. \cite{Vianna1} and \cite{Vianna2}).

In particular, the toric degenerations of $\cp^2$ are described by the Markov triples $(a,b,c)\in\Z^3$, $a\le b\le c$, satisfying to the Diophantine equation $a^2+b^2+c^2=3abc$ with the infinite number of solutions. Distinct Markov triples $(a,b,c)$ produce distinct monotone Lagrangian tori $L(a^2,b^2,c^2)$ coming from the monotone Clifford tori in the weighted projective plane $\Pp(a^2,b^2,c^2)$. 

For other Del Pezzo surfaces the toric degenerations with Picard number equal to 1 are described in \cite{HaPr} while those of Picard number greater than 1 are constructed in \cite{CMG}. All these degenerations produce exotic monotone Lagrangian tori in the corresponding Del Pezzo surfaces. Furthermore, Theorem \ref{thm-main} implies that symplectic geometry of the resulting torus detects the Picard number of the degeneration.
\end{exa}

\begin{exa}
For $n>2$, we can get examples of exotic Lagrangian tori by taking the product of tori from \ref{exa-dim2} and $(\cp^1,L_1)$, where $L_1$ is the monotone fiber circle in $\cp^1\approx S^2$. 
Theorem \ref{thm-main} implies that if two toric degeneration tori $L,L'\subset X$ were non-equivalent in a Del Pezzo surface $X$ then they stay non-equivalent in $X\times S^2\times\dots\times S^2$ after taking the product with $L_1$. 
\end{exa}

By the Gromov compactness, the number of classes $\alpha\in H_1(L)$ with $\nd(\alpha)\neq 0$ is finite. 
The element
$$\sum\nd(\alpha)\alpha\in\Z[H_1(L)]$$
of the group ring of $H_1(L)$ is called the {\em Floer potential}.

\begin{rmk} \label{rmk-terminal}
If the central fiber of the degeneration has terminal singularities,
then the Floer potential of the Lagrangian torus $L$ is determined by \Cref{thm-main},
because the only non-trivial lattice points of the polytope are vertices and the respective
number of discs for them equals one.

In the non-terminal case, the Theorem reduces the determination of Floer potentials to a finite algebraic question.
\end{rmk}

\begin{rmk} \label{rmk-batyrev}
The idea to use toric degenerations to study mirror symmetry for Fano manifolds
first appeared in works of Batyrev with Ciocan-Fontanine, Kim and van Straten \cite{BCFKS1,BCFKS2}.
In the case of a toric degeneration with terminal Gorenstein singularities of a central fiber,
Batyrev conjectured in \cite{Ba04,BK13} 
a mirror duality between the Fano manifold and
the Laurent polynomial 
described in \Cref{rmk-terminal}.
Our result shows that Batyrev's polynomials are Floer potentials of the monotone Lagrangian tori.
A specific form of mirror duality as discussed in \cite{EHX,BCFKS1,BCFKS2,Ba04,Ga-phd,Ga-ECM,Ga-sigma,BK13}
also follows by applying an open-closed relation, discovered in the context of toric degenerations by Bondal and Galkin \cite{BoGa},
and proved for any monotone Lagrangian torus by Tonkonog \cite{Tonkonog} 
with the help of the symplectic field theory.

\end{rmk}

In the presence of non-terminal singularities 
the central fiber alone is not sufficient to completely determine the Floer potential.
\begin{exa}
Let $X_0 = C S_6 \subset \Pp^7$ be the anti-canonical cone over the Del Pezzo surface of degree six $S_6 \subset \Pp^6$.
Recall that $S_6$ is a smooth toric Fano surface associated with a regular hexagon.
It is the closure in $\Pp^6$ of the embedding of $(\C^\times)^2$ of the linear system of the Laurent polynomials whose Newton polygon is contained in the hexagon $\Delta_6\subset\R^2$ with vertices $\pm(1,0)$, $\pm(0,1)$ and $\pm(1,1)$.
Accordingly, $X_0\subset\Pp^7$ is the toric 3-fold given by the pyramid of height 1 over $\Delta_6$ (treated as the Newton polytope for the linear system on $(\C^\times)^3$). 
Thus the polytope $\Delta^*$ for the toric 3-fold $X_0$ is a pyramid of height 2 over a hexagonal base with a single lattice non-vertex point.
In other words, for any $\qg$-smoothing $Y$ of $X_0$, \Cref{thm-main} determines all but one coefficients of the Floer potential
$$ W_{L_Y} = z + (x+y+x/y+y/x+1/x+1/y+c)/z,\ c\in\Z$$
of the monotone Lagrangian torus $L_Y$.

As it was noticed by Altmann \cite{Alt},
there are two distinct smoothings of $X_0$:
it can be smoothed to the product $\Pp^1 \times \Pp^1 \times \Pp^1$, or to 
the complete flag variety $F \subset \Pp^2 \times \Pp^2$.
Both smoothings are $\Q$-Gorenstein, but the coefficient $c$ is different: we have $c=3$ in the former case, and $c=2$ in the latter case, see \cite{Ga-sigma}.
the ambient spaces are not homotopically equivalent and have different Betti numbers.
\end{exa}

\begin{rmk}
The so-called {\em Minkowski
ansatz \cite{Ga-ECM,Ga-sigma}} assigns the binomial coefficient $\binom{l}{k}$
to the $k$th lattice point on the edge of the integer length $l$ of a lattice polytope.
One can extend the argument of the Theorem \ref{thm-main} to show that these binomial coefficients are equal
to the respective numbers of holomorphic disks. In this case the respective local model becomes the product
of $(n-2)$-dimensional complex torus with the two-dimensional $A_{l-1}$-singularity. $\C^2/(\Z_l)$.
\end{rmk}

Finally, we mention a recent high-dimensional application of Theorem \ref{thm-main} in \cite{BGM-sympl}.
We refer the reader to loc. cit. for details
and a slightly simplified adaptation of our method for this particular case.
\begin{exa} \label{exa-bgm}
Let $\mathcal{N}_g$ be
the moduli space of rank two stable bundles
with fixed determinant of odd degree over a Riemann surface of genus $g\geq 2$.
As a complex variety, it is a Fano manifold of complex dimension $3g-3$ with $H^2(\mathcal{N}_g) = \Z \Theta$ and
$c_1(\mathcal{N}_g) = 2 \Theta$.
As a symplectic manifold, it is monotone, and its symplectomorphism class depends only on $g$, see e.g. \cite{SiTi}.
For any trivalent connected graph $\gamma$ with $b_1(\gamma)=g$, Manon \cite{Manon} constructed a
$\qg$-degeneration of $\mathcal{N}_g$ to
the toric Fano variety $\Pp_\gamma$
of complex dimension $3g-3$, given by the spherical triangle inequalities associated to $\gamma$.
The monotone (Clifford) fiber torus of  $\Pp_\gamma$ thus yields a Lagrangian monotone torus $L_\gamma$ in $\mathcal{N}_g$ defined by the graph $\gamma$.

The combinatorial non-abelian Torelli theorem \cite{BGM-torelli} asserts that $\gamma$
is determined by the toric variety $\Pp_\gamma$, or equivalently by either its moment polytope $\Delta_\gamma$ or its dual fan polytope $\Delta_\gamma^\vee$.
If $\gamma$ is a loopless graph then $\Pp_\gamma$ has terminal singularities,
so \Cref{thm-main} determines the Floer potential by \Cref{rmk-terminal}.
In fact, the Floer potential for $(\mathcal{N}_g,L_\gamma)$ equals to graph potentials $W_\gamma$ introduced in \cite{BGM-tqft},
and thus $(\mathcal{N}_g,L_\gamma)$ and $W_\gamma$ are mirror dual.
In particular, different 3-valent graphs $\gamma$ produce different monotone Lagrangian tori in $\mathcal{N}_g$.

Another curious aspect of this example is that the total number of holomorphic disks of Maslov index 2 for $L_\gamma$ is $8g-8$, and thus depends only on the genus of the graph $\gamma$.
According to \cite{BGM-sympl}, it is the smallest possible number for any monotone Lagrangian torus in $\mathcal{N}_g$ .
\end{exa}

\noindent{\bf Acknowledgement} We thank Jo\'e Brendel and Felix Schlenk for their helpful feedback and corrections.

\bibliographystyle{plain}
\end{document}